\documentclass[11pt, a4paper]{article}

\usepackage{amsfonts}
\usepackage{amssymb,amsthm, amsmath}
\usepackage{bbm}

\input arrow.tex

\newtheorem{theorem}{Theorem}
\newtheorem{definition}[theorem]{Definition}
\newtheorem{proposition}[theorem]{Proposition}
\newtheorem{corollary}[theorem]{Corollary}
\newtheorem{lemma}[theorem]{Lemma}

\theoremstyle{remark}

\newtheorem{example}[theorem]{Example}

\newtheorem{remark}[theorem]{Remark}

\title{Simplicial complexes and minimal free resolution of monomial algebras}
\author{Ignacio Ojeda\thanks{During the preparation of this paper, both authors were partially supported by Ministerio de Educaci\'on y Ciencia (Spain), project MTM2007-65638.}\\ Universidad de Extremadura \\
Departamento de Matem\'{a}ticas\\
e-mail: ojedamc@unex.es \and
A. Vigneron-Tenorio$^*$\\
Universidad de C\'{a}diz\\
Departamento de Matem\'{a}ticas\\
e-mail: alberto.vigneron@uca.es}

\begin{document}

\date{\today}
\maketitle

\begin{abstract}
This paper is concerned with the combinatorial description of the graded minimal free resolution of certain monomial algebras which includes toric rings. Concretely, we explicitly describe how the graded minimal free resolution of those algebras is related to the combinatorics of some simplicial complexes. Our description may be interpreted as an algorithmic procedure to partially compute this resolution.


\smallskip
{\small \emph{MSC-class:} 16W50 (Primary) 13D02, 13F55 (Secondary).}
\end{abstract}

\section*{Introduction}

Let $I$ be an ideal in a polynomial ring $R$ over a field $\mathbbmss{k}.$ The $\mathbbmss{k}-$algebra $R/I$ is said to be monomial, if the algebraic set $\mathcal{V}(I)$ is parameterized by monomials.

Let $R/I$ be a monomial algebra. Since monomial algebras are semigroup algebras, one can consider a semigroup $S$ to study $R/I.$ This approach makes it possible to define a particular $S-$grading on the monomial algebra $R/I$ which allows to define the $S-$graded minimal free resolution of $R/I$ as $R-$module, under some reasonable hypothesis on $S$ (see Section \ref{Sect Prelim}). This graded minimal free resolution of $R/I$ has been explored by many authors with remarkable success (see e.g. \cite{Vigneron} and the references therein).

The study of the graded minimal free resolution of the monomial algebra $R/I$ from a semigroup viewpoint facilitates the use of methods based on the knowledge of the combinatorics of the semigroup. This paper is focused on this direction.

In this paper, we consider the simplicial complexes introduced by S.~Elihaou in his PhD Thesis \cite{Eliahou} and we show how their reduced $j-$th homology vector spaces over $\mathbbmss{k}$ are related with the $j-$th module of syzygies appearing in the graded minimal free resolution of certain monomial algebras (see Corollary \ref{Cor Iso}). Of course, this is not a very surprising theoretical result. A similar one was given by E.~Briales et al. in \cite{Collectanea}, although they used another different simplicial complexes. In fact, we prove that the reduced homology of both simplicial complexes are isomorphic (Theorem \ref{Th IsoHom}) and then, we use the results in \cite{Collectanea} to reach our Corollary \ref{Cor Iso}. Therefore, in this part, our main contribution should be regarded as showing the utility of the Elihaou's simplicial complexes for studying monomial algebras.

It is convenient to note here that, in some cases, Elihaou's simplicial complexes have a better behavior than the other ones and vice versa; for instance, when the minimal syzygies are concentrated in small $S-$degrees, Elihaou's simplicial complexes seems to be the right choice. This is the case of the monomial algebras $R/I$ such that $I$ generated by its indispensable binomials (see, e.g. \cite{Charalambous07, OjVi2}) which is of special interest in Algebraic Statistics and includes generic lattice ideals (\cite{Peeva}) and Lawrence type semigroup ideals (\cite{Sturmfels95, PisVi}).

The second part of the paper (Section \ref{Sect CSC}) is devoted to the explicit computation of minimal systems of generators of the $j-$th module of syzygies $N_j$ of a monomial algebra with associated semigroup $S.$ The main problem one encounters in the known combinatorial algorithms consists in determine \emph{a priori} the $S-$degrees in which the minimal generators of the syzygies are concentrated. This problem is not solved yet. Moreover, in this case, the information provided by the simplicial complexes used by E.~Briales et al. is clearly insufficient (see \cite{Collectanea, Campillo}). Therefore, we propose a different approach (Theorem \ref{Th DivComb0}, Proposition \ref{Prop DivComb1} and Corollary \ref{main_theorem}): we fix any $S-$degree $m$ and compute a subset of minimal generators of $N_j$ by only constructing one Elihaou's simplicial complex associated to $m$ and choosing suitable bases for some $\mathbbmss{k}-$vector spaces. This result does not solve the general problem, but opens the door to new perspectives in the combinatorial description of monomial (semigroup) algebras. In particular, an algorithmic procedure to compute a chain complex of free $R-$modules contained in minimal free resolution of a monomial algebra $R/I$ is described.

\section{Preliminaries}\label{Sect Prelim}

Let $S$ denote a commutative semigroup with zero element $0 \in S.$ Let $G(S)$ be a commutative group with a semigroup homomorphism $\iota : S \rightarrow G(S)$ such that every homomorphism from $S$ to a group factors in a unique way through $\iota.$ The commutative group $G(S)$ exists and is unique up to isomorphism, it is called the associated commutative group of $S.$ Further, $G(S)$ is finitely generated when $S$ is. The map $\iota$ is injective if, and only if, $S$ is cancellative, that is to say, if $m + n = m + n',\ m, n$ and $n' \in S,$ implies $n = n',$ in this case, $G(S)$ is the smallest group containing $S.$

For the purpose of this paper, we will assume that $S$ is combinatorially finite, i.e., there are only finitely many ways to write $m \in S$ as a sum $m = m_1 + \ldots + m_q,$ with $m_i \in S \setminus \{0\}.$ Equivalently, $S$ is combinatorially finite if, and only if, $S \cap (-S) = \{0\}$ (see Proposition 1.1 in \cite{BCMP}). Notice that this property guarantees that $m' \prec_S m \Longleftrightarrow m-m' \in S$ is a well defined partial order on $S.$

From now on, $S$ will denote a finitely generated, combinatorially finite, cancellative and commutative semigroup. We write $\mathbbmss{k}[S]$ for the $\mathbbmss{k}-$vector space $$\mathbbmss{k}[S] = \bigoplus_{m \in S} \mathbbmss{k} \chi^m$$ endowed with a multiplication which is $\mathbbmss{k}-$linear and such that $\chi^m \cdot \chi^n := \chi^{m+n},\ m$ and $n \in S.$ Thus $\mathbbmss{k}[S]$ has a natural $\mathbbmss{k}-$algebra structure and we will refer to it as the semigroup algebra of $S.$

In addition, we will fix a system of nonzero generators $n_1, \ldots, n_r$ for $S.$ Thus, $\mathbbmss{k}[S]$ may be regarded as the monomial $\mathbbmss{k}-$algebra generated by $\chi^{n_1}, \ldots, \chi^{n_r}.$

Moreover, this choice of generators induces a natural $S-$grading on $R = \mathbbmss{k}[x_1, \ldots, x_r],$ by assigning weight $n_i$ to $x_i,\ i = 1, \ldots, r,$ that is to say, $$R = \bigoplus_{m \in S} R_m,$$ where $R_m$ is the vector subspace of $R$ generated by all the monomials $\mathbf{x}^\alpha := x_1^{a_1} \cdots x_r^{a_r}$ with $\sum_{i=1}^r a_i n_i = m$ and $\alpha = (a_1, \ldots, a_n).$ Since $S$ is combinatorially finite, the vector spaces $R_m$ are finite dimensional (see Proposition 1.2 in \cite{BCMP}). We will denote by $\mathfrak{m}$ the irrelevant ideal of $R,$ that is to say, $\mathfrak{m} = \bigoplus_{m \in S \setminus \{0\}} R_m = (x_1, \ldots, x_r).$

\subsection{Minimal resolution}

The surjective $\mathbbmss{k}-$algebra morphism $$\varphi_0 : R \longrightarrow \mathbbmss{k}[S];\ x_i \longmapsto \chi^{n_i}$$ is $S-$graded, thus, the ideal $I_S := \ker(\varphi_0)$ is a $S-$homogeneous ideal called the ideal of $S.$ Notice that $I_S$ is a toric ideal (in the sense of \cite{Sturmfels95} chapter 4) generated by $$\Big\{\mathbf{x}^\alpha - \mathbf{x}^\beta : \sum_{i=1}^r a_i n_i = \sum_{i=1}^r b_i n_i \in S \Big \}.$$

Now, by using the $S-$graded Nakayama's lemma recursively (see Proposition 1.4 in \cite{BCMP}), we may construct $S-$graded $\mathbbmss{k}-$algebra homomorphism $$\varphi_{j+1} : R^{s_{j+1}} \longrightarrow R^{s_j},$$ corresponding to a choice of a minimal set of $S-$homogeneous generators for each module of syzygies $N_j := \ker(\varphi_j),$ notice that $N_0 = I_S.$ Thus, we obtain a minimal free $S-$graded resolution for the $R-$module $\mathbbmss{k}[S]$ of
type
$$
\ldots \longrightarrow R^{s_{j+1}} \stackrel{\varphi_{j+1}}{\longrightarrow} R^{s_j} \longrightarrow \ldots \longrightarrow R^{s_2} \stackrel{\varphi_2}{\longrightarrow} R^{s_1} \stackrel{\varphi_1}{\longrightarrow} R
\stackrel{\varphi_{0}}{\longrightarrow} \mathbbmss{k}[S] \longrightarrow 0,
$$
where $s_{j+1} := \sum_{m \in S} \mathrm{dim}_\mathbbmss{k} V_j(m),$ with $V_j(m) := (N_j)_m/(\mathfrak{m} N_j)_m,$ is the so-called $(j+1)-$th Betti number. Observe that the dimension of $V_j(m)$ is the number of generators of degree $m$ in a minimal system of generators of the $j-$th module of syzygies $N_j$ (i.e. the multigraded Betti number $s_{j,m}$), so, by the Noetherian property of $R,\ s_{j+1}$ is finite. Moreover, the Auslander-Buchsbaum's formula assures that $s_j = 0$ for $j > p = r - \mathrm{depth}_R \mathbbmss{k}[S]$ and $s_p \neq 0.$ (cf. Theorem 1.3.3 in \cite{BH}).

\subsection{Simplicial homology}\label{SSect SH}

Let $K$ be a finite simplicial complex on $[n] := \{1, \ldots, n\}.$ For each integer $i,$ let $\mathcal{F}_i(K)$ be the set of $i-$dimensional faces of $K,$ and let $\mathbbmss{k}^{\mathcal{F}_i(K)}$ be a $\mathbbmss{k}-$vector space whose basis element $\mathbf{e}_F$ correspond to $i-$faces $F \in \mathcal{F}_i(K).$

The reduced chain complex of $K$ over $\mathbbmss{k}$ is the complex $\widetilde{C}_\bullet(K):$
$$
0 \rightarrow \mathbbmss{k}^{\mathcal{F}_{n-1}(K)} \stackrel{\partial_{n-1}}{\longrightarrow} \ldots \longrightarrow \mathbbmss{k}^{\mathcal{F}_{i}(K)} \stackrel{\partial_i}{\longrightarrow} \mathbbmss{k}^{\mathcal{F}_{i-1}(K)} \longrightarrow \ldots \stackrel{\partial_0}{\longrightarrow} \mathbbmss{k}^{\mathcal{F}_{-1}(K)} \to 0
$$
The boundary maps $\partial_i$ are defined by setting $\mathrm{sing}(j,F) = (-1)^{r-1}$ if $j$ is the $r-$th element of the set $F \subseteq [n],$ written in increasing order, and
$$
\partial_i(\mathbf{e}_F) = \sum_{\substack{j \in F \\ \# F = i+1}} \mathrm{sing}(j,F) \mathbf{e}_{F \setminus j}.
$$

For each integer $i,$ the $\mathbbmss{k}-$vector space
$$
\widetilde{H}_i(K) = \ker(\partial_i)/\mathrm{im}(\partial_{i+1})
$$
in homological degree $i$ is the $i-$th reduced homology of $K.$ Elements of $\widetilde{Z}_i(K) := \ker(\partial_i)$ are called $i-$cycles and elements of $\widetilde{B}_i(K) := \mathrm{im}(\partial_{i+1})$ are called $i-$boundaries.

\section{Simplicial complexes and minimal syzygies}

In this section, we will consider two different simplicial complexes associated with $S$ and we will compare their homologies. The first simplicial complex was introduced by S. Eliahou in \cite{Eliahou} and the second one is used in \cite{Collectanea} to describe the minimal free resolution of $\mathbbmss{k}[S].$

For any $m \in S,$ let $C_m = \{\mathbf{x}^\alpha = x_1^{a_1} \cdots x_r^{a_r} \mid \sum_{i=1}^r a_i n_i = m\}$ and define the abstract simplicial complex on the vertex set $C_m,$ $$\nabla_m = \{ F \subseteq C_m \mid \gcd(F) \neq 1 \},$$ where $\gcd(F)$ denotes the greatest common divisor of the monomials in $F.$ Notice that $\nabla_m$ has finitely many vertices because $S$ is combinatorially finite.

For any $m \in S,$  we consider the abstract simplicial complex on the vertex set $[r],$ $$\Delta_m = \{F \subseteq [r] \mid m - n_F \in S\},$$ where $n_F = \sum_{i \in F} n_i.$

Now, we are going to compare $\widetilde{H}_\bullet(\nabla_m)$ with $\widetilde{H}_\bullet(\Delta_m).$ To facilitate our work, we recall the so-called ``Nerve Lemma".

\begin{definition}
A cover of a simplicial complex $K$ is a family of subcomplexes $\mathcal{K} = \{K_\alpha \mid \alpha \in A\}$ with $K = \bigcup_{\alpha \in A} K_\alpha.$

We say that the cover $\mathcal{K}$ satisfies the Leray property if each non-empty finite intersection $K_{\alpha_1} \cap \ldots \cap K_{\alpha_q}$ is acyclic.
\end{definition}

\begin{definition}
Let $\mathcal{K} = \{K_\alpha \mid \alpha \in A\}$ be a cover of a simplicial complex $K.$ The nerve of $\mathcal{K},$ denoted by $N_\mathcal{K},$ is the simplicial complex having vertices $A$ and with $\{\alpha_1, \ldots, \alpha_q\}$ being a simplex if $\bigcap_{i=1}^q K_{\alpha_i} \neq \varnothing.$
\end{definition}

\medskip
\noindent\textbf{Nerve Lemma.} \emph{Assume that $\mathcal{K} = \{K_\alpha \mid \alpha \in A\}$ is a cover of a simplicial complex $K.$ If $\mathcal{K}$ satisfies the Leray property, then $$ H_j (N_\mathcal{K}) \cong H_j(K),$$ for all $j \geq 0.$}

\medskip
\begin{proof}
See Theorem 7.26 in \cite{Rotman}.
\end{proof}

\begin{theorem}\label{Th IsoHom}
$\widetilde{H}_j (\nabla_m) \cong  \widetilde{H}_j(\Delta_m),$ for all $j \geq 0$ and $m \in S.$
\end{theorem}

\begin{proof}
For each $\mathbf{x}^\alpha \in C_m,$ define the simplicial complex $K_\alpha = \mathcal{P}(\mathrm{supp}(\mathbf{x}^\alpha)),$ that is to say, the full subcomplex of $\Delta_m$ with vertices $\mathrm{supp}(\mathbf{x}^\alpha).$ Set $\mathcal{K}^m =\{K_\alpha :\mathbf{x}^\alpha \in C_m\}.$

On the one hand, we have that $F \in \Delta_m,$ i.e. , $m - n_F \in S$ if, and only, if, there exists $\mathbf{x}^\alpha \in C_m$ with $\mathrm{supp}(\mathbf{x}^\alpha) \supseteq F,$ therefore, $\mathcal{K}^m$ is a cover of $\Delta_m.$

Moreover, $\bigcap_{i=1}^q K_{\alpha_i} \neq \varnothing$ if, and only if, $\mathrm{gcd}(\mathbf{x}^{\alpha_1}, \ldots, \mathbf{x}^{\alpha_q}) \neq 1,$ so $\nabla_m$ is the nerve of $\mathcal{K}^m.$ Finally, since the cover $\mathcal{K}^m$ of $\Delta_m$ satisfies the Leray property, because $\bigcap_{i=1}^q K_{\alpha_i} \neq \varnothing$ is a full simplex, by the Nerve Lemma, we may conclude the existence of the desired isomorphism.
\end{proof}

The above theorem has been proved independently by H.~Charalambous and A.~Thoma (see Theorem 3.2 in \cite{Charalambous08}).

\medskip
By Theorem 2.1 in \cite{Collectanea}, one has that $V_j(m) \cong \widetilde{H}_j (\Delta_m),$ for all $m \in S.$ Therefore, we have the following elementary consequence:

\begin{corollary}\label{Cor Iso}
$\widetilde{H}_j (\nabla_m) \cong V_j(m),$ for all $j \geq 0$ and $m \in S.$
\end{corollary}

Notice that the above corollary assures that the multigraded Betti number $s_{j+1, m}$ equals the rank of the $j-$reduced homology group $\widetilde{H}_j(\nabla_m)$ of the simplicial complex $\nabla_m,$ for every $m \in S.$

Furthermore, we emphasize that Corollaries 2.2 and 2.3 in \cite{Collectanea} may be written in terms of the complexes $\nabla_\bullet,$ by simply using Theorem \ref{Th IsoHom}. Obtaining by this way necessary and sufficient combinatorial conditions for $\mathbbmss{k}[S]$ to be Cohen-Macaulay or Gorenstein. Indeed, $\mathbbmss{k}[S]$ is Cohen-Macaulay if, and only if, $\widetilde{H}_{r-d}(\nabla_m) = 0,$ for every $m \in S,$ where $d = \mathrm{rank}(G(S)).$ In this case, the Cohen-Macaulay type of $\mathbbmss{k}[S]$ is $$s_{r-d} = \sum_{m \in S} \dim \widetilde{H}_{r-d-1}(\nabla_m).$$

\section{On the computation of $\widetilde{H}_{j}(\nabla_m).$}\label{Sect CHj}

One of the keys to our results in the next section consist in the assumption that we are able to compute (and fix) a particular basis for the $\mathbbmss{k}-$vector space $\widetilde{Z}_j(\nabla_m),$ for each $j \geq -1$ and $m \in S.$

To do this we consider the reduced chain complex $\widetilde{C}_\bullet(\nabla_m),\ m \in S,$ as defined in Subsection \ref{SSect SH} and order the faces according to a (fixed) criterion, e.g. by choosing a monomial term order $\prec$ on $R.$ Indeed, $\prec$ induces a well ordering on the $j-$dimensional faces: $F < F'$ if, and only if, $\gcd(F) \prec \gcd(F').$

Thus, by decreasingly ordering all the $j-$dimensional faces according to the chosen criterion a basis $\mathcal{B}_j = \{F^{(j)}_1, \ldots, F^{(j)}_{d_j}\}$ of $\mathbbmss{k}^{\mathcal{F}_j(\nabla_m)}$ is fixed, for each $j \geq 0$ and $m \in S.$

Let $A_j \in \mathbb{Z}^{d_{j-1} \times d_j}$ be the matrix of $\partial_j$ with respect to $\mathcal{B}_j$ and $\mathcal{B}_{j-1},\ j \geq 0.$ By performing Gaussian elimination on $A_j$ two invertible matrices $P_j$ and $Q_j$ are obtained such that $$P_j^{-1} A_j Q_j = \left(\begin{array}{c|c} I_{r_j} & 0 \\ \cline{1-2} 0 & 0 \end{array}\right) \in \mathbb{Z}^{d_{j-1} \times d_j},$$ where $I_{r_j}$ is the identity matrix of order $r_j = \mathrm{rank}(A_j),\ j \geq 0.$ Then, the first $r_j$ columns of $P_j$ are the coordinates with respect to $\mathcal{B}_{j-1}$ of a basis of $\widetilde{B}_{j-1}(\nabla_m) = \mathrm{im}(\partial_j)$ and the last $d_j - r_j$ columns of $Q_j$ are the coordinates with respect to $\mathcal{B}_j$ of a basis of $\widetilde{Z}_j(\nabla_m) = \ker(\partial_j),$ for each $j \geq 0$ and $m \in S.$ 

Now, since $\widetilde{B}_j(\nabla_m) \subseteq \widetilde{Z}_j(\nabla_m),$ by using the bases obtained above and elementary linear algebra, we can extend the basis of $\widetilde{B}_j(\nabla_m)$ to a basis of $\widetilde{Z}_j(\nabla_m),$ for each $j \geq 0$ and $m \in S.$

Therefore, we may construct a $\mathbbmss{k}-$basis
\begin{equation}\label{ecu basis}\big\{\widehat{\mathbf{h}}^{(j)}_1, \ldots, \widehat{\mathbf{h}}^{(j)}_{t'_j}, \widehat{\mathbf{b}}^{(j)}_1, \ldots, \widehat{\mathbf{b}}^{(j)}_{t''_j}\}\end{equation} of $\widetilde{Z}_j(\nabla_m)$ such that
\begin{itemize}
\item[(a)] $\widehat{\mathbf{h}}^{(j)}_i = \sum_{k=1}^{d_j} q^{(j)}_{ki} \partial_{j+1}\big(F^{(j+1)}_k\big),$ where $q^{(j)}_{ki}$ is the $(k,i)-$th entry of $Q_j,\ i = 1, \ldots, t'_j.$
\item[(b)] the classes of $\widehat{\mathbf{b}}^{(j)}_1, \ldots, \widehat{\mathbf{b}}^{(j)}_{t''_j}$ modulo $\widetilde{B}_j(\nabla_m)$ form a $\mathbbmss{k}-$basis of $\widetilde{H}_j(\nabla_m),$
\end{itemize}
for each $j \geq 0$ and $m \in S.$

\begin{remark}\label{Rem Z0}
Since $A_0 = (1\ 1\ \ldots\ 1) \in \mathbb{Z}^{1 \times d_0},$ we may assume that the corresponding basis for $\widetilde{Z}_0(\nabla_m)$ is $\big\{ \{\mathbf{x}^{\beta_1}\} - \{\mathbf{x}^\alpha\}, \ldots, \{\mathbf{x}^{\beta_{d_0}}\} - \{\mathbf{x}^\alpha\} \big\},$ with $\mathbf{x}^\alpha \succ \mathbf{x}^{\beta_1} \succ \ldots \succ \mathbf{x}^{\beta_{d_0}}.$ So, $\widehat{\mathbf{b}}^{(0)}_i = \{\mathbf{x}^{\beta_{k_i}}\} - \{\mathbf{x}^\alpha\}$ for some $k_i \in \{1, \ldots, d_0\}.$
\end{remark}

Notice that this general construction can be also applied to compute a $\mathbbmss{k}-$basis of $\widetilde{H}_j(\Delta_m).$ In any case, the computation of $\widetilde{H}_j(\nabla_m)$ and $\widetilde{H}_j(\Delta_m)$ is equally difficult (see \cite{Campillo} for a different approach on the computation of $\widetilde{H}_j(\Delta_m)$).

\section{Computing syzygies from combinatorics}\label{Sect CSC}

In this section, we will explicitly describe the isomorphisms whose existence we have proved in Corollary \ref{Cor Iso}.

We will start by giving an isomorphism $\widetilde{H}_0(\nabla_m) \stackrel{\sigma_0}{\cong} V_0(m).$ As the reader can note, the construction of $\sigma_0$ follows from the definition of $\widetilde{H}_0(\nabla_m) = \widetilde{Z}_0(\nabla_m)/\widetilde{B}_0(\nabla_m)$ and $V_0(m) = (N_0)_m/(\mathfrak{m} N_0)_m.$ However, we will give the construction by taking in mind the general case in order to introduce the notation of this section.

\medskip
First of all, consider the $\mathbbmss{k}-$linear map
\begin{equation}\label{eq2.1}
\displaystyle{\psi_0:\mathbbmss{k}^{\mathcal{F}_0(\nabla_m)}} \longrightarrow R;\
\{\mathbf{x}^\alpha\} \longmapsto \mathbf{x}^\alpha.
\end{equation}
This map induces an isomorphism from $\widetilde{Z}_0(\nabla_m)$ to $(N_0)_m.$ More precisely,
$$
\widetilde{Z}_0(\nabla_m) \longrightarrow (N_0)_m;\ \widehat{b} := \{\mathbf{x}^\alpha\} - \{\mathbf{x}^\beta\} \longmapsto b := \mathbf{x}^\alpha - \mathbf{x}^\beta,
$$
recall that $\widetilde{Z}_0(\nabla_m)$ is generated by $\{\mathbf{x}^\alpha\} - \{\mathbf{x}^\beta\},$ with $\mathbf{x}^\alpha$ and $\mathbf{x}^\beta \in C_m$ (see Remark \ref{Rem Z0}), and $(N_0)_m$ is generated by pure difference binomials of $S-$degree equals $m.$

Therefore, we have a surjective map $\overline{\psi}_0$ given by the composition
$$
\widetilde{Z}_0(\nabla_m) \longrightarrow  (N_0)_m \longrightarrow V_0(m) = (N_0)_m/(\mathfrak{m} N_0)_m.
$$

\begin{lemma}\label{Lema IsoH0}
$\widetilde{B}_0(\nabla_m) \subseteq \ker \overline{\psi}_0.$
\end{lemma}

\begin{proof}
Since $\widehat{f} \in \widetilde{B}_0(\nabla_m) = \mathrm{im}(\partial_1),$ there exist $\{\mathbf{x}^{\alpha_j}, \mathbf{x}^{\beta_j}\} \in \nabla_m$ and $\mu_j \in \mathbbmss{k},$ such that $$\partial_1 \left(\sum_j \mu_j \{\mathbf{x}^{\alpha_j}, \mathbf{x}^{\beta_j}\} \right) =  \sum_j \mu_j \left(\{\mathbf{x}^{\beta_j}\} - \{\mathbf{x}^{\alpha_j}\} \right) = \widehat{f}.$$ So,
$$
f = \overline{\psi}_0 (\widehat{f})= \sum_j \mu_j \left(\mathbf{x}^{\beta_j} - \mathbf{x}^{\alpha_j} \right) = \sum_j
\mu_j \mathbf{x}^{\gamma_j} \left(\underbrace{\mathbf{x}^{\beta'_j} -
\mathbf{x}^{\alpha'_j}}_{\in N_0} \right),
$$
with $\mathbf{x}^{\gamma_j} = \mathrm{gcd}(\mathbf{x}^{\alpha_j}, \mathbf{x}^{\beta_j}), \mathbf{x}^{\alpha'_j} = \mathbf{x}^{\alpha_j}/\mathbf{x}^{\gamma_j}$ and $\mathbf{x}^{\beta'_j} = \mathbf{x}^{\beta_j}/\mathbf{x}^{\gamma_j},$ moreover, $\mathbf{x}^{\gamma_j} \neq 1,$ because $\{ \mathbf{x}^{\alpha_j}, \mathbf{x}^{\beta_j}\}$ is an edge of $\nabla_m.$ Thus, 
$f \in (\mathfrak{m} N_0)_m$ as claimed.
\end{proof}

Therefore, by Lemma \ref{Lema IsoH0}, $\overline{\psi}_0$ factorizes canonically through $\widetilde{H}_0(\nabla_m):$

$$
\sarrowlength=.5\harrowlength \commdiag{ \widetilde{Z}_0(\nabla_m) & \mapright\lft{\overline{\psi}_0} & V_0(m) \cr &
\arrow(1,-1)\rt{\pi} \quad \arrow(1,1)\rt{\sigma_0}\cr & \widetilde{H}_0(\nabla_m) \cr}
$$
Notice, that $\sigma_0$ is an isomorphism because it is surjective and, by Corollary \ref{Cor Iso}, $\dim \widetilde{H}_0(\nabla_m) = \dim V_0(m).$

\medskip
Now we will show a combinatorial method to compute some minimal binomial generators of $I_S$ from a given binomial in $I_S.$ But first, we will introduce an important property of the complexes $\nabla_m$ which claims that $\nabla_{m'}$ can be easily computed from $\nabla_m,$ for every $m' \prec_S m,$ i.e., if $m-m' \in S.$

\begin{lemma}\label{Lema S-less}
Let $m$ and $m' \in S.$ If $m' \prec_S m,$ then
$$
\nabla_{m'} \cong \{ F \in \nabla_m \mid \mathbf{x}^\beta\ \mbox{properly divides}\ \gcd(F) \},
$$
for any (fixed) monomial $\mathbf{x}^\beta \in C_{m-m'}.$
\end{lemma}

\begin{proof}
Let $\mathbf{x}^\beta$ a monomial in $C_{m-m'}.$ If $F' = \{\mathbf{x}^{\alpha'_1}, \ldots, \mathbf{x}^{\alpha'_t}\} \in \nabla_{m'},$ then $F = \{\mathbf{x}^{\alpha'_1 + \beta}, \ldots, \mathbf{x}^{\alpha'_t + \beta}\} \in \nabla_m.$ Conversely, consider $F = \{\mathbf{x}^{\alpha_1}, \ldots,$ $\mathbf{x}^{\alpha_t}\} \in \nabla_m$ such that $\mathbf{x}^\beta$ divides $\gcd(F)$ and $\gcd(F) \neq \mathbf{x}^\beta.$ Since $\mathbf{x}^\beta$ divides $\mathbf{x}^{\alpha_i},\ i = 1, \ldots, t,$ and $\gcd(F) \neq \mathbf{x}^\beta,$ we conclude that $F' = \{\mathbf{x}^{\alpha_1 - \beta}, \ldots, \mathbf{x}^{\alpha_t - \beta}\}$ is a face of $\nabla_{m'}.$
\end{proof}

\begin{theorem}\label{Th DivComb0}
Let $m \in S$ and let $\nabla_m$ be given. For each $\mathbf{x}^\alpha - \mathbf{x}^\beta \in (I_S)_m,$ it can be computed a unique subset $\mathcal{B} = \{b_1, \ldots, b_t\}$ of a minimal system of binomial generators of $I_S$ and unique $f_1, \ldots, f_t \in R,$ such that
\begin{itemize}
\item[(a)] $\mathbf{x}^\alpha -\mathbf{x}^\beta = \sum_{i=1}^t f_i b_i,$
\item[(b)] $\gcd(\mathbf{x}^\alpha, \mathbf{x}^\beta)$ divides $f_i,\ i = 1, \ldots, t.$
\end{itemize}
\end{theorem}

\begin{proof}
We divide the proof in two steps.

\medskip
\noindent STEP 1.- Write $$\mathbf{x}^\alpha - \mathbf{x}^\beta = \mathbf{x}^\gamma (\mathbf{x}^{\alpha'} - \mathbf{x}^{\beta'}),$$ where $\mathbf{x}^\gamma = \gcd(\mathbf{x}^\alpha, \mathbf{x}^\beta),\ \mathbf{x}^{\alpha'} = \mathbf{x}^\alpha/\mathbf{x}^\gamma$ and $\mathbf{x}^{\beta'} = \mathbf{x}^\beta/\mathbf{x}^\gamma.$ Notice that $\{\mathbf{x}^{\alpha}\}$ and $\{\mathbf{x}^{\beta}\}$ are adjacent in $\nabla_m$ when $\mathbf{x}^\gamma \neq 1,$ and that $\{\mathbf{x}^{\alpha'}\}$ and $\{\mathbf{x}^{\beta'}\}$ are never adjacent in $\nabla_{m'},$ where $m'$ is the $S-$degree of $\mathbf{x}^{\alpha'}$ (and $\mathbf{x}^{\beta'},$ of course). Moreover, $m' \prec_S m,$ when $\mathbf{x}^\gamma \neq 1.$ In this case, we consider the simplicial complex $\nabla_{m'}$ (computed from $\nabla_m$ by using Lemma \ref{Lema S-less}) and the binomial $\mathbf{x}^{\alpha'} - \mathbf{x}^{\beta'} \in (I_S)_{m'}.$

\medskip
\noindent STEP 2.- For simplicity, by Step 1, we may assume that $\{\mathbf{x}^{\alpha}\}$ and $\{\mathbf{x}^{\beta}\}$ are not adjacent.

Let $\{\widehat{h}_1, \ldots, \widehat{h}_{t'},\ \widehat{b}_1, \ldots, \widehat{b}_{t''} \}$ be a $\mathbbmss{k}-$basis of $\widetilde{Z}_0(\nabla_m)$ constructed as in Section \ref{Sect CHj}. Then, $\widehat{h}_j = \sum_{k=1}^{d_1} q^{(0)}_{kj} \partial_1(F_k^{(1)}) \in \widetilde{B}_0(\nabla_m),$ for every $j,$ where $\mathcal{F}_1(\nabla_m) = \{F^{(1)}_1, \ldots, F^{(1)}_{d_1}\},$ and the classes of $\widehat{b}_1, \ldots, \widehat{b}_{t''}$ modulo $\widetilde{B}_0(\nabla_m)$ form a basis of $\widetilde{H}_0(\nabla_m).$

Set $\widehat{g}_k = \partial_1(F_k^{(1)}) = \{\mathbf{x}^{\alpha_{k}}\} - \{\mathbf{x}^{\beta_{k}}\}.$ Then,
$$
\{\mathbf{x}^\alpha\} - \{\mathbf{x}^\beta\} = \sum_i \lambda_i \widehat{b}_i + \sum_j \mu_j \widehat{h}_j =
\sum_i \lambda_i \widehat{b}_i + \sum_k \big(\sum_j \mu_j  q^{(0)}_{kj}\big) \widehat{g}_k
,
$$
for some $\lambda_i$ and $\mu_j \in \mathbbmss{k}.$

By Remark \ref{Rem Z0}, we have that $\widehat{b}_i,$ is a pure difference of vertices in $\nabla_m$ for every $i = 1, \ldots, t''.$ Therefore,
\begin{equation}\label{eq3.3}
\mathbf{x}^\alpha - \mathbf{x}^\beta =\psi_0 (\{\mathbf{x}^\alpha\} - \{\mathbf{x}^\beta\}) =  \sum_i
\lambda_i b_i + \sum_k \nu_k g_k,
\end{equation}
where the ${b_i}$'s are binomials in $R$ and $\nu_k = \sum_j \mu_j  q^{(0)}_{kj} \in \mathbbmss{k},\ k = 1, \ldots, d_1.$

If $\nu_k = 0,$ for every $k,$ we are done. Otherwise, we repeat this procedure (starting from Step 1) for each $g_k = \mathbf{x}^{\alpha_{k}} - \mathbf{x}^{\beta_{k}}$ with $\nu_k \neq 0.$ Since $\gcd(\mathbf{x}^{\alpha_{k}},\mathbf{x}^{\beta_{k}}) \neq 1,$ the $S-$degree of the binomial produced in Step 1 will be strictly lesser than the degree of $g_k,$ so, we may guarantee that this process ends in finitely many iterations\footnote{It is convenient to recall that, by Lemma \ref{Lema S-less}, we do not need to compute the new complexes $\nabla_{m'_k},$ for any $k.$}.

Finally, notice that we have considered a particular basis for each $\widetilde{Z}_0(\nabla_\bullet)$ appearing in. Thus, our computation depends on the choice of these bases. However, no other choice has been made. Thus, assuming fixed basis for each $\widetilde{Z}_0(\nabla_\bullet)$ (see Section \ref{Sect CHj}), we may guarantee that $b_i$'s and $f_i$'s are uniquely obtained.
\end{proof}

\begin{remark}\label{Rem PostTh}
Observe that the proof of Theorem \ref{Th DivComb0} may be considered as an algorithm which effectively computes
a subset of a minimal set of binomial generators of $I_S$ starting from any binomial of $S-$degree $m.$ 

Further, if the sets $\mathcal{B}$'s corresponding to two different binomials in $(I_S)_m$ have an element of the same $S-$degree, then they have the same elements of such $S-$degree (because the bases of $\widetilde{Z}_0(\nabla_\bullet)$ are chosen fixed), that is to say, both sets are subsets of the same minimal system of binomial generators of $I_S.$

Thus, by considering the union of all sets $\mathcal{B}$'s corresponding to each binomial in $(I_S)_m$ a subset of a minimal system of binomial generators of $I_S$ is obtained.
\end{remark}

\medskip
We show with an example how to compute some minimal binomial generators by using the above theorem.

\begin{example}\label{ejemplo1}
Let $S \subset \mathbb{Z}^2$ be the semigroup generated by the columns of the following matrix
$$
\left(
\begin {array}{cccc} 4 & 5 & 7 & 8 \smallskip \\  1 & 1 & 1 & 1 \end{array}
\right).
$$
The binomial $g = {x_2}^2{x_3}^6-{x_1}^3{x_4}^5$ is
clearly in $I_S \subset \mathbbmss{k }[x_1, \ldots, x_4]$ and its
$S-$degree is $m = (52,8) = 2\, (5,1) + 6\, (7,1) \in S.$

The vertex set\footnote{One can compute this set by solving diophantine equations (see \cite{Nsol}).} of $\nabla_m$ is
\begin{align*}
C_m = \Big\{ & {x_2}^2{x_3}^6,\ {x_2}^3{x_3}^3{x_4}^2,\ {x_2}^4{x_4}^4,\ {x_1}{x_2}{x_3}^5{x_4},\ {x_1}{x_2}^2{x_3}^2{x_4}^3,\\ & {x_1}^2{x_3}^4{x_4}^2,\ {x_1}^2{x_2}{x_3}{x_4}^4,\ {x_1}^3{x_4}^5
\Big\}
\end{align*}

The simplicial complex $\nabla_m$ is clearly connected, so $$\Big\{ \{\mathbf{x}^\alpha\} - \{ {x_2}^3{x_3}^3{x_4}^2\} \mid \mathbf{x}^\alpha \in C_m \setminus \{{x_2}^3{x_3}^3{x_4}^2\} \Big\}$$ is a $\mathbbmss{k}-$basis of $\widetilde{B}_0(\nabla_m) = \widetilde{Z}_0(\nabla_m)$ and
\begin{align*}
\widehat{g} & = \{ {x_2}^2{x_3}^6 \} - \{ {x_1}^3{x_4}^5 \}\\ & =  \left( \{ {x_2}^2{x_3}^6 \} - \{ {x_2}^3{x_3}^3{x_4}^2 \} \right) - \left( \{{x_1}^3{x_4}^5\} - \{ {x_2}^3{x_3}^3{x_4}^2 \}\right)
\end{align*}

Since ${x_2}^2{x_3}^6 - {x_2}^3{x_3}^3{x_4}^2 = {x_2}^2{x_3}^3 \left({x_3}^3 - {x_2}{x_4}^2\right)$ and ${x_1}^3{x_4}^5 - {x_2}^3{x_3}^3{x_4}^2$ $= {x_4}^2 \left({x_1}^3{x_4}^3 - {x_2}^3{x_3}^3 \right).$ The binomials that we are interested in now are $$g_1 = {x_3}^3 - {x_2}{x_4}^2\quad \text{and}\quad g_2 = {x_1}^3{x_4}^3 - {x_2}^3{x_3}^3$$ which have $S-$degree $m_1 = (21,3)$ and $m_2 = (36,6) \in S,$ respectively.

By Lemma \ref{Lema S-less}, $C_{m_1} = \{{x_3}^3 ,\ {x_2}{x_4}^2\}$ and $\nabla_{m_1}$ is disconnected. Thus, we
conclude that $b_1 = g_1$ is minimal binomial generator of $I_S.$ On the other hand, by Lemma \ref{Lema S-less} again, we have that
\begin{align*}
C_{m_2} = \{ & {x_2}^3{x_3}^3,\ {x_2}^4{x_4}^2,\ {x_1}{x_2}^2{x_3}^2{x_4},\ {x_1}^2{x_3}^4, {x_1}^2{x_2}{x_3}{x_4}^{2},\ {x_1}^3{x_4}^3\}
\end{align*}
and that $\nabla_{m_2}$ is connected. A $\mathbbmss{k}-$basis of
$\widetilde{B}_0(\nabla_{m_2}) = \widetilde{Z}_0(\nabla_{m_2})$ is
$$\Big\{ \{\mathbf{x}^\alpha\} - \{{x_1}{x_2}^2{x_3}^2{x_4}\} \mid \mathbf{x}^\alpha \in C_{m_2} \setminus \{{x_1} {x_2}^2{x_3}^2{x_4}\} \Big\}$$
and
\begin{align*}
\widehat{g}_2 & = \{{x_2}^3{x_3}^3 \}  - \{{x_1}^3{x_4}^3\} = \\ & =  \Big( \{{x_2}^3{x_3}^3 \} -
 \{ {x_1}{x_2}^2{x_3}^2{x_4}\} \Big) - \Big(\{{x_1}^3{x_4}^3\} - \{{x_1}{x_2}^2{x_3}^2{x_4}\}\Big).
\end{align*}

Since ${x_2}^3{x_3}^3 - {x_1}{x_2}^2{x_3}^2{x_4} = {x_2}^2 {x_3}^2 \big({x_2}{x_3} - {x_1}{x_4}\big)$ and ${x_1}^3{x_4}^3 - {x_1}{x_2}^2{x_3}^2 \cdot $ $\cdot {x_4} = x_1 x_4 \big({x_1}^2{x_4}^2-{x_2}^2{x_3}^2\big)$ we have to consider now the  binomials $g_{21}={x_2}{x_3} - {x_1}{x_4}$ and $g_{22} = {x_1}^2{x_4}^2-{x_2}^2{x_3}^2.$ By using Lemma \ref{Lema S-less} in order to compute the corresponding simplicial complexes, it is easy to see that $b_2 = g_{21} \in (I_S)_{(12,2)}$ is a minimal binomial generator and that $g_{22} = {x_1}^2{x_4}^2-{x_2}^2{x_3}^2 \in (I_S)_{(24,4)}$ is not, because $$C_{m_3} = \{{x_1}^2{x_4}^2, {x_1}{x_2}{x_3}{x_4}, {x_2}^2{x_3}^2\},$$ with $m_3 = (24,2),$ and $\nabla_{m_3}$ is connected. A $\mathbbmss{k}-$basis of $\widetilde{B}_0(\nabla_{m_2}) = \widetilde{Z}_0(\nabla_{m_2})$ is $\big\{ \{{x_1}^2{x_4}^2\}-\{{x_1}{x_2}{x_3}{x_4}\}, \{{x_2}^2{x_3}^2\} - \{{x_1}{x_2}{x_3}{x_4}\} \big\}$ and then
\begin{align*}
\widehat{g}_{22} & = \{{x_1}^2{x_4}^2\} - \{{x_2}^2{x_3}^2\} \\
& = \Big(\{{x_1}^2{x_4}^2\}-\{{x_1}{x_2}{x_3}{x_4}\}\Big) - \Big(\{{x_2}^2{x_3}^2\} - \{{x_1}{x_2}{x_3}{x_4}\}\Big)
\end{align*}
Thus, $g_{22} = {x_1}{x_4}({x_1}{x_4} - {x_2}{x_3}) - {x_2}{x_3}({x_2}{x_3} - {x_1}{x_4}) = (-{x_1}{x_4}-{x_2}{x_3}) ({x_2}{x_3} - {x_1}{x_4})$

Summarizing, we have obtained two minimal binomial generators,
$b_1 = {x_3}^3 - x_2{x_4}^2$ and $b_2 = {x_2}{x_3}-{x_1}{x_4}$ of $I_S$ and two polynomials $f_1 = {x_2}^2{x_3}^3,$ and $f_2 = x_4^2 \big( {x_2}^2{x_3}^2 + {x_1}{x_2}{x_3}{x_4}$ $+$ ${x_1}^2 {x_4}^2\big)$ such that
$$
{x_2}^2{x_3}^6-{x_1}^3{x_4}^5= f_1 b_1 + f_2 b_2.
$$
\end{example}

\subsection{First syzygies of semigroup ideals}

\begin{remark}
For the sake of simplicity in the notation, we will assume that we have obtained a whole minimal system of binomial generators of $I_S$ by using Theorem \ref{Th DivComb0}:
$$
\mathcal{B} := \{b_1, \ldots, b_{s_1}\}
$$
Although this will be not truly necessary for our purpose. In practice, we will only need to know the subset of minimal binomial generators of $I_S$ obtained by applying Theorem \ref{Th DivComb0} to each $\mathbf{x}^\alpha - \mathbf{x}^\beta \in (I_S)_m$ with $\gcd(\mathbf{x}^\alpha, \mathbf{x}^\beta) \neq 1$ (see Remark \ref{Rem PostTh}).
\end{remark}


\medskip
We define the $\mathbbmss{k}-$linear map $$\psi_1:\displaystyle{\mathbbmss{k}^{\mathcal{F}_1(\nabla_m)}} \longrightarrow R^{s_1};\ \{\mathbf{x}^\alpha, \mathbf{x}^\beta\} \longmapsto \mathbf{f} := \left(\begin{array}{c} f_1 \\ \vdots \\ f_{s_1}\end{array}\right)$$ where $f_i \neq 0$ is given by Theorem \ref{Th DivComb0} from $\mathbf{x}^\alpha - \mathbf{x}^\beta \in (I_S)_m.$ Since we are working with fixed bases for $\widetilde{Z}_j(\nabla_\bullet),\ j \geq 0,$ we may assure that $\psi_1$ is well defined.

\begin{example}
For instance, in Example \ref{ejemplo1}, we have obtained that
$$\psi_1 \left(\left\{ {x_2}^2{x_3}^6 ,\ {x_1}^3{x_4}^5 \right\}\right)= \left(\begin{array}{c} {x_2}^2{x_3}^3 \\ {x_2}^2{x_3}^2{x_4}^2 + {x_1}{x_2}{x_3}{x_4}^2 + {x_1}^2 {x_4}^4\\ 0 \\ \vdots \\ 0 \end{array}\right) \in R^{s_1}.$$
\end{example}

\begin{remark}
In the following, we will write $[f_1, \ldots, f_s] \in R^s$ instead of to use column-vector notation.
\end{remark}

\begin{lemma}
The map $\psi_1$ makes commutative the following diagram
\begin{equation}\label{diag1}
\commdiag{
\displaystyle{\mathbbmss{k}^{\mathcal{F}_1(\nabla_m)}} & \mapright\lft{\psi_1} & R^{s_1} \cr \mapdown\lft{\partial_1} & & \mapdown\rt{\varphi_1} \cr \displaystyle{\mathbbmss{k}^{\mathcal{F}_0(\nabla_m)}} & \mapright\lft{\psi_0} &
R, \cr}\end{equation}
where the bottom row map is defined as in (\ref{eq2.1}) and $\varphi_1(\mathbf{e}_i) = b_i,\ i = 1, \ldots, s_1.$
\end{lemma}

\begin{proof}
Consider $\widehat{\mathbf{f}} = \sum_{j=1}^{d_1} \lambda_j \{\mathbf{x}^{\alpha_j}, \mathbf{x}^{\beta_j}\} \in \mathbbmss{k}^{\mathcal{F}_1(\nabla_m)}.$ Then, \begin{align*}
\varphi_1 \circ \psi_1(\widehat{\mathbf{f}}) & = \varphi_1 \Big( \sum_{j=1}^{d_1} \lambda_j [f_{1j}, \ldots, f_{s_1 j}] \Big) = \sum_{j=1}^{d_1} \lambda_j \sum_{i=1}^{s_1} f_{ij} b_i\\ & = \sum_{j=1}^{d_1} \lambda_j (\mathbf{x}^{\alpha_j} - \mathbf{x}^{\beta_j}) = \psi_0 \Big( \sum_{j=1}^{d_1} \lambda_j \big(\{\mathbf{x}^{\alpha_j}\} - \{\mathbf{x}^{\beta_j}\}\big) \Big)\\ & = \psi_0 \circ \partial_1 \Big( \sum_{j=1}^{d_1} \lambda_j \big\{\mathbf{x}^{\alpha_j}, \mathbf{x}^{\beta_j}\big\} \Big) = \psi_0 \circ \partial_1(\widehat{\mathbf{f}}).
\end{align*}
\end{proof}

Furthermore, one can see that $\psi_1$ sends $1-$cycles to $1-$syzygies. Indeed, if $z = \sum_i \lambda_i
\{\mathbf{x}^{\alpha_i},\mathbf{x}^{\alpha_{i+1}}\} \in \widetilde{Z}_1(\nabla_m),$ then $\partial_1(z) = 0.$ Thus, if $\psi_1(z) = \sum_i \lambda_i [f_{i1},\ldots, f_{is_1}],$ we have that
\begin{equation}\label{syz1}
0 = \varphi_1\Bigl(\sum_i \lambda_i [f_{i1},\ldots, f_{is_1}] \Bigr) = \sum_i \lambda_i \sum_j f_{ij} b_j = \sum_j \Bigl(\sum_i \lambda_i f_{ij} \Bigr) b_j.
\end{equation}
Thus, $ [\sum_i \lambda_i f_{i1},\ldots, \sum_i \lambda_i f_{is_1}]$ is a syzygy, as claimed. The converse is also true in the following sense:

\begin{lemma}\label{Lema ClaveH1.1}
The map $\psi_1:\widetilde{Z}_1(\nabla_m) \longrightarrow  (N_1)_m$ is surjective.
\end{lemma}

\begin{proof}
Let $\mathbf{f} = [f_1, \ldots, f_{s_1}] \in (N_1)_m.$ Thus, if $b_i = \mathbf{x}^{\alpha_i} - \mathbf{x}^{\beta_i}$ and $f_i = \sum_j \lambda_{ij} \mathbf{x}^{\gamma_{ij}},\ i = 1, \ldots, s_1,$ it follows that
\begin{align*}
0 & = \sum_i f_i b_i = \sum_i \Bigl(\sum_j \lambda_{ij} \mathbf{x}^{\gamma_{ij}}\Bigr) \Bigl(\mathbf{x}^{\alpha_i} - \mathbf{x}^{\beta_j} \Bigr) \\ & = \sum_{i, j} \lambda_{ij} \Bigl(\mathbf{x}^{\gamma_{ij}} \mathbf{x}^{\alpha_i} - \mathbf{x}^{\gamma_{ij}} \mathbf{x}^{\beta_j} \Bigr).
\end{align*}
By taking $\widehat{\mathbf{f}} = \displaystyle{\sum_{i,j}} \lambda_{ij} \{\mathbf{x}^{\gamma_{ij} + \beta_i}, \mathbf{x}^{\gamma_{ij} + \alpha_i} \}\in \widetilde{Z}_1(\nabla_m),$ we conclude that $\mathbf{f} = \psi_1 (\widehat{\mathbf{f}}).$
\end{proof}

\medskip
Thus, we have a surjective map $\overline{\psi}_1$ which is nothing but the composition
$$
\widetilde{Z}_1(\nabla_m) \longrightarrow  (N_1)_m \longrightarrow V_1(m) = (N_1)_m/(\mathfrak{m} N_1)_m.
$$

\begin{lemma}\label{Lema ClaveH1.2}
$\widetilde{B}_1(\nabla_m) \subseteq \ker \overline{\psi}_1.$
\end{lemma}

\begin{proof}
Since $\widetilde{B}_1(\nabla_m) = \mathrm{im}(\partial_2)$ and $\partial_2$ is $\mathbbmss{k}-$linear, it suffices to prove that $\partial_2(F) \in \ker \overline{\psi}_1$ for any $2-$dimensional face $F$ of $\nabla_m.$

Let $F = \{\mathbf{x}^{\alpha_1}, \mathbf{x}^{\alpha_2}, \mathbf{x}^{\alpha_3}\}$ be a $2-$dimensional face of $\nabla_m.$ Then $\mathbf{x}^\gamma = \mathrm{gcd}(F) \neq 1.$ Thus, by Theorem \ref{Th DivComb0}, there exist $[f_{i1}, \ldots, f_{is_1}] \in (N_1)_m,\ i = 1, 2, 3,$ such that $\mathbf{x}^\gamma$ divides $f_{ij}$ and
\begin{align*}
\overline{\psi}_1(\partial_2(F)) & = \overline{\psi}_1\left(\{\mathbf{x}^{\alpha_2}, \mathbf{x}^{\alpha_3}\} - \{\mathbf{x}^{\alpha_1}, \mathbf{x}^{\alpha_3}\} + \{\mathbf{x}^{\alpha_1}, \mathbf{x}^{\alpha_2}\}\right) \\
& = [f_{11}, \ldots, f_{1s_1}] - [f_{21}, \ldots, f_{2s_1}] + [f_{31}, \ldots, f_{3s_1}] \\ & = \mathbf{x}^\gamma \Bigl([f'_{11}, \ldots, f'_{1s_1}] - [f'_{21}, \ldots, f'_{2s_1}] + [f'_{31}, \ldots, f'_{3s_1}]\Bigr) \\ & = \mathbf{x}^\gamma [f'_{11} - f'_{21} + f'_{31}, \ldots, f'_{1s_1} - f'_{2s_1} + f'_{3s_1}]
\end{align*}
Therefore, $\overline{\psi}_1(\partial_2(F)) \in (\mathfrak{m} N_1)_m,$ as claimed.
\end{proof}

By Lemma \ref{Lema ClaveH1.2}, $\overline{\psi}_1$ factorizes canonically through $\widetilde{H}_1(\nabla_m):$ $$\sarrowlength=.5\harrowlength \commdiag{ \widetilde{Z}_1(\nabla_m) & \mapright\lft{\overline{\psi}_1} & V_1(m) \cr & \arrow(1,-1)\rt{\pi} \quad \arrow(1,1)\rt{\sigma_1}\cr & \widetilde{H}_1(\nabla_m) \cr}$$ As before, $\sigma_1$ is an isomorphism because it is surjective and, by Corollary \ref{Cor Iso}, $\dim \widetilde{H}_1(\nabla_m) = \dim V_1(m).$

\begin{proposition}\label{Prop DivComb1}
Let $m \in S$ and let $\nabla_m$ be given. For each $\mathbf{g} := [g_1, \ldots,$ $g_{s_1}] \in (N_1)_m,$  it can be computed a unique subset $\mathcal{B} = \{\mathbf{b}^{(1)}_1, \ldots, \mathbf{b}^{(1)}_t\}$ of a minimal system of generators of the first module of syzygies of $\mathbbmss{k}[S]$ and unique $f_1, \ldots, f_t \in R$ such that
\begin{itemize}
\item[(a)]$\mathbf{g} = \sum_{j=1}^t f_j \mathbf{b}^{(1)}_j,$
\item[(b)] $\gcd(g_1, \ldots, g_{s_1})$ divides $f_j,\  j = 1, \ldots, t.$
\end{itemize}
\end{proposition}

\begin{proof}
Write $$\mathbf{g} = h\, \mathbf{g}',$$ where $h := \gcd(g_1, \ldots, g_{s_1}).$ Notice that when $h \neq 1,$ the $S-$degree of $\mathbf{g}'$ is strictly lesser than the $S-$degree of $\mathbf{g}.$ In this case, we consider the simplicial complex $\nabla_{m'},$ where $m'$ is the $S-$degree of $\mathbf{g}'$ (recall that $\nabla_{m'}$ can be computed by using Lemma \ref{Lema S-less}) and $\mathbf{g}' = [g'_1, \ldots, g'_{s_1}] \in (N_1)_{m'}.$

\medskip
For simplicity, we assume that $\gcd(g_1, \ldots, g_{s_1}) = 1,$ i.e. $\mathbf{g} = \mathbf{g}'.$

\medskip
Let $\{\widehat{\mathbf{h}}_1, \ldots, \widehat{\mathbf{h}}_{t'},\ \widehat{\mathbf{b}}_1, \ldots, \widehat{\mathbf{b}}_{t''}\}$ be a $\mathbbmss{k}-$basis of $\widetilde{Z}_1(\nabla_m)$ constructed as in Section \ref{Sect CHj}. Then $\widehat{\mathbf{h}}_j = \sum_{k=1}^{d_2} q^{(1)}_{kj} \partial_2\big(F_k^{(2)}\big) \in \widetilde{B}_1(\nabla_m),$ for every $j,$ and the classes of $\widehat{\mathbf{b}}_1, \ldots, \widehat{\mathbf{b}}_{t''}$ form a basis of $\widetilde{H}_1(\nabla_m).$

So, if $\widehat{\mathbf{g}} \in \widetilde{Z}_1(\nabla_m)$ is such that $ \mathbf{g}=\psi_1(\widehat{\mathbf{g}}) \in
(N_1)_m$ (see Lemma \ref{Lema ClaveH1.1} and its proof), then
$$
\widehat{\mathbf{g}} = \sum_i \lambda_i \widehat{\mathbf{b}}_i + \sum_j \mu_j \widehat{\mathbf{h}}_j = \sum_i \lambda_i \widehat{\mathbf{b}}_i + \sum_k \big(\sum_j \mu_j q^{(1)}_{kj} \big) \widehat{\mathbf{g}}_k
$$
for some $\lambda_i$ and $\mu_j \in \mathbbmss{k},$ with $\widehat{\mathbf{g}}_k = \partial_2\big(F_k^{(2)}\big).$ Therefore,
\begin{equation*}
\mathbf{g} = \psi_1 (\widehat{\mathbf{g}})= \sum_i \lambda_i \mathbf{b}_i + \sum_k \nu_k \mathbf{g}_k,
\end{equation*}
where $\nu_k = \sum_j \mu_j  q^{(1)}_{kj} \in \mathbbmss{k},\ k = 1, \ldots, d_2.$

If $\nu_k = 0,$ for every $k,$ we are done. Otherwise, we repeat this procedure for $\mathbf{g}_k$ with $\nu_k \neq 0.$ Since $\gcd(F_k^{(2)}) \neq 1$ divides $\gcd(\mathbf{g}_k),$ we may assure that this process ends for $S-$degree reasons.

The uniqueness follows from the same argument as in the proof of Theorem \ref{Th DivComb0}.
\end{proof}

\begin{remark}
Similarly to the case of Theorem \ref{Th DivComb0}, we may assure that the subsets $\mathcal{B}$'s are contained in the same minimal system of generators of the first module of syzygies of $\mathbbmss{k}[S]$ (see Remark \ref{Rem PostTh}).
\end{remark}

\medskip
Let us illustrate the above theorem with an example.

\begin{example}\label{ejemplo2}
Let $\mathcal{B}= \{  x_2x_3-x_1x_4,\ x_3^3-x_2x_4^2,\ x_1x_3^2-x_2^2x_4,\ x_2^3-x_1^2x_3 \}$ be a minimal generating set of the semigroup ideal $I_S$ of example \ref{ejemplo1}. Consider $m = (45,7) \in S$ and
$$
\mathbf{g} = [{x_2}{x_3}^4 + {x_1}{x_2}{x_4}^3,-{x_2}^2{x_{{3}}}^2-{x_1}^2{x_4}^2, {x_2}{x_3}^2{x_4} + {x_1}{x_3}{x_4}^2,0] \in (N_1)_m.
$$
Let $\widehat{\mathbf{g}} \in \widetilde{Z}_1(\nabla_{(45, 7)}),$ such that $\mathbf{g}= \psi_1(\widehat{\mathbf{g}}),$ be defined as in the proof of  Lemma \ref{Lema ClaveH1.1}:
\begin{align*}
\widehat{\mathbf{g}} = & \big\{ {x_2}^2{x_3}^5,\ {x_1}{x_2}{x_3}^4{x_4} \big\} + \big\{{x_1}{x_2}^2{x_3}{x_4}^3,\ {x_1}^2{x_2}{x_4}^4  \big\} \\ & - \big\{{x_2}^2{x_3}^2{x_4}^2,\ {x_2}{x_3}^5 \big\} - \big\{ {x_1}^2{x_3}^3{x_4}^{ 2},\ {x_1}^2{x_2}{x_4}^4 \big\}+ \ldots
\end{align*}
As in the proof of Proposition \ref{Prop DivComb1}, we fix a particular basis of $\widetilde{Z}_1(\nabla_{(45,7)})$ in such way we may write
\begin{align*}
\widehat{\mathbf{g}} = & \Big(\big\{{x_2}^2{x_3}^5,{x_1}{x_2}{x_3}^4{x_4}\big\} - \big\{ {x_2}^2{x_3}^2{x_4}^2,{x_2}{x_3}^5\big\}\\ & + \big\{ {x_1}{x_2}^2{x_3}{x_4}^2 ,{x_1}^2{x_3}^3{x_4} \big\} \Big) + \Big( \big\{ {x_1}{x_2}^2{x_3}{x_4}^3, {x_1}^2{x_2}{x_4}^4 \big\}\\ & - \big\{ {x_1}^2{x_3}^3{x_4}^{2}, {x_1}^2{x_2}{x_4}^4 \big\} + \big\{ {x_1}{x_2}{x_3}^4{x_4}, {x_2}^3{x_3}^2{x_4}^2 \big\}\Big) + \ldots
\end{align*}
The image by $\psi_1$ of the first parenthesis is
$$
{x_3}[-{x_2}{x_3}^3, {x_2}^2{x_3},-{x_1}{x_4}^2,0] \in N_1.
$$
Then by taking $\mathbf{g}_1 = [-{x_2}{x_3}^3, {x_2}^2{x_3}, -{x_1}{x_4}^2,0] \in (N_1)_{(38,6)}$ we may repeat the above process again and so. By proceeding similarly with all the other parenthesis, we finally get:
\begin{align*}
\mathbf{g} = &{x_2}{x_3}^2[-{x_3}^2,{x_2},-{x_4},0] + {x_1}{x_4}^2[-{x_2}{x_4},{x_1},-{x_3},0]
\end{align*}
Thus, we have obtained two minimal syzygies in $S-$degrees $(25,4)$ and $(26,4),$ respectively.
\end{example}

\subsection{$i-$syzygies of semigroup ideal}

Let $m \in S$ and let $\nabla_m$ be given. Let us suppose that
\begin{itemize}
\item We are able to compute a $\mathbbmss{k}-$linear map $\psi_{i-1} :  \displaystyle{\mathbbmss{k}^{\mathcal{F}_{i-1}(\nabla_m)}} \longrightarrow R^{s_{i-1}}$ such that $\psi_i : \widetilde{Z}_{i-1}(\nabla_m) \to (N_{i-1})_m$ is a well defined surjective $\mathbbmss{k}-$linear map and $\widetilde{B}_{i-1}(\nabla_m) \subseteq \ker \overline{\psi}_{i-1},$ where $\overline{\psi}_{i-1}$ is the composition $\widetilde{Z}_{i-1}(\nabla_m) \longrightarrow  (N_{i-1})_m \longrightarrow V_{i-1}(m) = (N_{i-1})_m/(\mathfrak{m} N_{i-1})_m.$
\item For each $\mathbf{g} := [g_1, \ldots, g_{s_{i-1}}] \in (N_{i-1})_m,$ we are able to compute a unique subset $\{\mathbf{b}^{(i-1)}_1, \ldots, \mathbf{b}^{(i-1)}_t\}$ of a minimal system of generators of the $(i-1)-$th mo\-du\-le of syzygies of $\mathbbmss{k}[S]$ and unique $f_1, \ldots, f_t \in R$ such that $\mathbf{g} = \sum_{j=1}^t f_j \mathbf{b}^{(i-1)}_j$ and $\gcd(g_1, \ldots, g_{s_{i-1}})$ divides $f_j,\ j = 1, \ldots, t.$
\end{itemize}

Similarly to the former cases, we assume that a set, $\{\mathbf{b}^{(i-1)}_1, \ldots, \mathbf{b}^{(i-1)}_{s_i}\},$ of minimal generators of the $i-$th module of syzygies of $\mathbbmss{k}[S]$ is obtained from the above hypothetical computation (as before, a subset of the system will be enough for our needs).

Then, we may define the new $\mathbbmss{k}-$linear map
\begin{equation}\label{diag2}
\psi_i:\displaystyle{\mathbbmss{k}^{\mathcal{F}_i(\nabla_m)}} \longrightarrow R^{s_i};\ F := \{\mathbf{x}^{\alpha_0}, \ldots, \mathbf{x}^{\alpha_i} \} \longmapsto \mathbf{f} := [f_1, \ldots, f_{s_i}],
\end{equation}
where $f_i \neq 0$ is given by the above hypothetic computation from $\psi_{i-1}\circ\partial_i(F) \in (N_{i-1})_m.$

Thus, the map $\psi_i$ makes commutative the following diagram
\begin{equation}\label{diag3}
\commdiag{\displaystyle{\mathbbmss{k}^{\mathcal{F}_i(\nabla_m)}} & \mapright\lft{\psi_i} & R^{s_i} \cr \mapdown\lft{\partial_i} & & \mapdown\rt{\varphi_i} \cr \displaystyle{\mathbbmss{k}^{\mathcal{F}_{i-1}(\nabla_m)}} & \mapright\lft{\psi_{i-1}} & R^{s_{i-1}}, \cr}
\end{equation}
where $\varphi_1(\mathbf{e}_k) = b_k,\ k = 1, \ldots, s_{i-1}.$

%

As in (\ref{syz1}), it is easy to see that $\psi_i$ sends $i-$cycles to $i-$syzygies. Besides $\psi_i :
\widetilde{Z_i}(\nabla_m) \to (N_i)_m$ is surjective. Indeed, given $\mathbf{g} := [g_1, \ldots, g_{s_{i}}] \in (N_{i})_m$ with $g_j=\sum_k\lambda_{jk}\mathbf{x}^{\delta_{jk}}\in R,$ one obtains that$$0=\sum_j g_j \mathbf{b}^{(i-1)}_j= \sum_j \sum_k \lambda_{jk}\mathbf{x}^{\delta_{jk}}\mathbf{b}^{(i-1)}_j,$$ where $\mathbf{x}^{\delta _{jk}}\mathbf{b}^{(i-1)}_j\in (\mathfrak{m} N_{i-1})_m.$ Now, let us consider $$\widehat{\mathbf{x}^{\delta_{jk}}\mathbf{b}^{(i-1)}_j}= {\psi_{i-1}}^{-1} (\mathbf{x}^{\delta _{jk}}\mathbf{b}^{(i-1)}_j)\in \widetilde{B}_{i-1}(\nabla_m).$$ This cycle can be write as $\widehat{\mathbf{x}^{\delta_{jk}}\mathbf{b}^{(i-1)}_j}=\sum_l \mu_{jkl} \partial_{i} (\underbrace{F^{(i)}_l}_{\in \mathcal{F}_i(\nabla_m)}).$ Therefore, we conclude that
$$
\widehat{\mathbf{g}}=\displaystyle{\sum_j \sum_k \lambda_{jk} \sum_l \mu_{jkl} F^{(i)}_l} \in \widetilde{Z}_i (\nabla_m)
$$
satisfies $\mathbf{g}= \psi_i (\widehat{\mathbf{g}}),$ as claimed.

Finally, one can prove that $\widetilde{B}_i(\nabla_m) \subseteq \ker \overline{\psi}_i$ with the same arguments as in Lemma \ref{Lema ClaveH1.2}. So, we have surjection $\overline{\psi}_i :\widetilde{Z}_i(\nabla_m) \longrightarrow  (N_i)_m \longrightarrow V_i(m),$ which factorizes canonically through $\widetilde{H}_i(\nabla_m):$
$$
\sarrowlength=.5\harrowlength \commdiag{ \widetilde{Z}_i(\nabla_m) & \mapright\lft{\overline{\psi}_i} & V_i(m) \cr & \arrow(1,-1)\rt{\pi} \quad \arrow(1,1)\rt{\sigma_i}\cr & \widetilde{H}_i(\nabla_m) \cr}
$$
This defines an isomorphism $\sigma_i$ as desired.

\medskip
The last ingredient in our construction is the following result which guarantees that we will be able to define $\psi_{i+1}$ in similar terms as we assumed to be possible for $\psi_i.$ This will complete our main objective: to give an explicit description of the isomorphism in Corollary \ref{Cor Iso}.

\begin{corollary}\label{main_theorem}
Let $m \in S$ and let $\nabla_m$ be given. For each $\mathbf{g} := [g_1, \ldots, g_{s_{i}}] \in (N_{i})_m,$ it can be computed a unique subset $\{\mathbf{b}^{(i)}_1, \ldots, \mathbf{b}^{(i)}_t\}$ of a minimal system of generators of the $i-$th mo\-du\-le of syzygies of $\mathbbmss{k}[S]$ and unique $f_1, \ldots, f_t \in R$ such that
\begin{itemize}
\item[(a)] $\mathbf{g} = \sum_{j=1}^t f_j \mathbf{b}^{(i)}_j,$
\item[(b)] $\gcd(g_1, \ldots, g_{s_{i}})$ divides $f_j,\ j = 1, \ldots, t.$
\end{itemize}
\end{corollary}

\begin{proof}
This proof is the natural generalization of the proofs of Theorem
\ref{Th DivComb0} and Proposition \ref{Prop DivComb1}.

\medskip
Write $\mathbf{g} = h\, \mathbf{g}',$ where $h := \gcd(g_1, \ldots, g_{s_i}).$ The $S-$degree$(\mathbf{g}')\prec_S S-$degree$(\mathbf{g}),$ when $h \neq 1.$ For the sake of notation, we suppose that $\gcd(\mathbf{g}) = 1,$ i.e., $\mathbf{g} = \mathbf{g}'.$

Let $\{\widehat{\mathbf{h}}_1, \ldots, \widehat{\mathbf{h}}_{t'},\ \widehat{\mathbf{b}}_1, \ldots, \widehat{\mathbf{b}}_{t''}\}$ be a $\mathbbmss{k}-$basis of $\widetilde{Z}_i(\nabla_m)$ constructed as in Section \ref{Sect CHj}. Then $\widehat{\mathbf{h}}_j = \sum_{k=1}^{d_i} q^{(i)}_{kj} \partial_{i+1}\big(F_k^{(i+1)}\big) \in \widetilde{B}_i(\nabla_m),$ for every $j,$ and the classes of $\widehat{\mathbf{b}}_1, \ldots, \widehat{\mathbf{b}}_{t''}$ form a basis of $\widetilde{H}_i(\nabla_m).$

Now, we compute $\widehat{\mathbf{g}}\in \widetilde{Z}_i(\nabla_m)$ such that $\mathbf{g}= \psi_i( \widehat{\mathbf{g}})$ with $\widehat{\mathbf{g}} = \sum_i \lambda_i \widehat{\mathbf{b}}_i + \sum_k \big(\sum_j \mu_j q^{(i)}_{kj}\big) \widehat{\mathbf{g}}_k,$ for some $\lambda_i$ and $\mu_j \in \mathbbmss{k},$ with $\widehat{\mathbf{g}}_k = \partial_{i+1}(F_k^{(i+1)}).$ Therefore $ \mathbf{g} = \sum_i \lambda_i \mathbf{b}_i + \sum_k \nu_k \mathbf{g}_k,$ with $\nu_k = \sum_j \mu_j q^{(i)}_{kj}.$ If $\nu_k = 0,$ for every $k,$ we are done. Otherwise, we repeat this procedure for $\mathbf{g}_k$ with $\mu_k \neq 0.$ 
This process ends for $S-$degree reasons.

The uniqueness follows from the same argument as in the proof of Theorem \ref{Th DivComb0}.
\end{proof}

By the above corollary, we can conclude that, starting from any $i-$syzygy $\mathbf{g}$ of $I_S,$ our combinatorial algorithm computes a subset $\mathcal{B}'$ of a minimal generating set $\mathcal{B}$ of $N_i$ and the polynomial coefficients of $\mathbf{g}$ with respect $\mathcal{B}'$ (and therefore with respect to $\mathcal{B})$ without knowing other $i-$syzygies. It is very important to note that $\mathbf{g}$ is not relevant by itself. Given its $S-$degree $m,$ we can effectively produce $i-$syzygies of $\mathbbmss{k}[S]$ in the $S-$degree $m$ and subsets of minimal generators of $N_j,\ j \le i.$ All the construction lies in the simplicial complex $\nabla_m.$

\medskip
Finally, let us see how our algorithm produces part of the minimal free resolution of a semigroup algebra $\mathbbmss{k}[S]$ starting from one $S-$de\-gree. In fact, in the next example, we will get the whole resolution.

\begin{example}
Let $S$ be the semigroup in example \ref{ejemplo1} and consider $m = (60,10) \in S.$ The set of vertices of the simplicial complex $\nabla_m$ is
\begin{align*}
C_m =  \{ & {x_2}^5{x_3}^5,\ {x_2}^6{x_3}^2{x_4}^2,\ {x_1}{x_2}^4{x_3}^4{x_4},\ {x_1}{x_2}^5{x_3}{x_4}^3,\\ & {x_1}^2{x_2}^2{x_3}^6,\ {x_1}^2{x_2}^3{x_3}^3{x_4}^2,\ {x_1}^2{x_2}^4{x_4}^4,\ {x_1}^3{x_2}{x_3}^5{x_4},\\ & {x_1}^3{x_2}^2{x_3}^2{x_4}^3,\  {x_1}^4{x_3}^4{x_4}^2,\ {x_1}^4{x_2}{x_3}{x_4}^4,{x_1}^5{x_4}^5 \}
\end{align*}

We are going to ``capture" syzygies of $\mathbbmss{k}[S]$ using the method described in this section. To do that we choose the following $3-$dimensional face of $\nabla_m$
\begin{align*}
F = & \big\{ \underbrace{{x_1}{x_2}^4{x_3}^4{x_4}}_A,\ \underbrace{{x_1}^2{x_2}^2{x_3}^6}_{B},\ \underbrace{{x_1}^2{x_2}^3{x_3}^3{x_4}^2}_C,\ \underbrace{{x_1}^3{x_2}{x_3}^5{x_4}}_D \big\}
\end{align*}

By Theorem \ref{Th DivComb0}, considering the $1-$dimensional faces of $F,$ we are able to construct four $\mathbbmss{k}-$linearly independent minimal binomial generators of $I_S,$ that is, $0-$syzygies of $\mathbbmss{k}[S]$
\begin{align*}
\{A,B\} \longrightarrow\ b_1 & = {x_2}^2{x_4}-{x_1}{x_3}^2\\
\{A,C\} \longrightarrow\ b_2 & = {x_2}{x_3} - {x_1}{x_4}\\
\{A,D\} \longrightarrow\ b_3 & = {x_2}^3-{x_1}{x_3}{x_4}\\
\{B,C\} \longrightarrow\ b_4 & = {x_3}^3-{x_2}{x_4}^2\\
\{B,D\} \longrightarrow\ b_2 & \\
\{C,D\} \longrightarrow\ b_1 &
\end{align*}
Recall that the obtained coefficients are also needed, although we do not write them here.

Now, by Proposition \ref{Prop DivComb1}, using the $2-$dimensional faces of $F$ and the non-written above coefficients,
we are able to produce four $\mathbbmss{k}-$linearly independent $1-$syzygies of $\mathbbmss{k}[S]:$
\begin{align*}
\{A,B,C\} \longrightarrow\ \mathbf{b}_1^{(1)} & = [x_3, -x_2 x_4, 0, x_1, 0, \ldots, 0] \in R^{s_1}\\
\{A,B,D\} \longrightarrow\ \mathbf{b}_2^{(1)} & = [x_2, x_1 x_3, -x_4, 0, 0, \ldots, 0] \in R^{s_1}\\
\{A,C,D\} \longrightarrow\ \mathbf{b}_3^{(1)} & = [x_1, x_2^2, -x_3, 0, 0, \ldots, 0] \in R^{s_1}\\
\{B,C,D\} \longrightarrow\ \mathbf{b}_4^{(1)} & = [x_4, -x_3^2, 0, x_2, 0, \ldots, 0] \in R^{s_1}
\end{align*}
Again, we do not write here the obtained coefficients, although we insist that they are necessary to go further. Notice that the coordinates of $\mathbf{b}_i^{(1)},\ i = 1, \ldots, 4,$ has been completed with zeroes, because \emph{a priori} we do not know whether the rank of $R^{s_1}$ is $4.$

Finally, by Corollary \ref{main_theorem}, using the $3-$dimensional face of $F$ and the non-written above coefficients, we get one $2-$syzygy of $\mathbbmss{k}[S]:$
$$
F = \{A,B,C,D\} \longrightarrow\ \mathbf{b}_1^{(2)} = [-x_2, x_3, -x_4, x_1, 0, \ldots, 0] \in R^{s_2}
$$

Therefore, we have obtained a chain complex of free $R-$modules
\begin{equation}\label{ecu cc}
0 \to R \to R^4 \to R^4 \to R \to R/J,
\end{equation}
where $J = (b_1, b_2, b_3, b_4) \subset R,$ which is a subcomplex of the minimal free resolution of $\mathbbmss{k}[S].$

In this case, it is not difficult to see that $I_S = J$ is a Gorenstein ideal of codimension $2$ and thus (\ref{ecu cc}) is its minimal free resolution. \end{example}

\medskip \noindent
\textbf{Acknowledgments.-}
We would like to thank Prof. Emilio Briales-Morales and Prof. Antonio Campillo-L\'opez for helpful comments and suggestions. We want also to thank the anonymous referee for his/her comments.

\end{document}